\documentclass[11pt]{article}
\usepackage{amssymb,amsfonts,amsmath,latexsym,epsf,tikz,url}

\newtheorem{theorem}{Theorem}[section]

\newtheorem{corollary}[theorem]{Corollary}
\newtheorem{lemma}[theorem]{Lemma}
\newtheorem{remark}[theorem]{Remark}

\newcommand{\proof}{\noindent{\bf Proof.\ }}
\newcommand{\qed}{\hfill $\square$\medskip}

\begin{document}

\title{Distinguishing number and distinguishing index of certain graphs}

\author{
Saeid Alikhani  $^{}$\footnote{Corresponding author}
\and
Samaneh Soltani
}

\date{\today}

\maketitle

\begin{center}
Department of Mathematics, Yazd University, 89195-741, Yazd, Iran\\
{\tt alikhani@yazd.ac.ir, s.soltani1979@gmail.com}
\end{center}


\begin{abstract}
The distinguishing number (index) $D(G)$ ($D'(G)$) of a graph $G$ is the least integer $d$
such that $G$ has an vertex labeling (edge labeling)  with $d$ labels  that is preserved only by a trivial
automorphism. In this paper we compute these two parameters for some specific graphs.   Also we study the distinguishing number and the distinguishing index of corona product of two graphs.  
\end{abstract}


\medskip
\noindent{\bf AMS Subj.\ Class.:} 05C15, 05E18

\section{Introduction}

Let $G=(V,E)$ be a simple graph. We use the standard graph notation (\cite{Sandi}). In particular, $Aut(G)$ denotes the automorphism group of $G$.  
A labeling of $G$, $\phi : V \rightarrow \{1, 2, \ldots , r\}$, is said to be $r$-distinguishing, 
if no non-trivial  automorphism of $G$ preserves all of the vertex labels.
The point of the labels on the vertices is to destroy the symmetries of the
graph, that is, to make the automorphism group of the labeled graph trivial.
Formally, $\phi$ is $r$-distinguishing if for every non-trivial $\sigma \in Aut(G)$, there
exists $x$ in $V = V (G)$ such that $\phi(x) \neq \phi(x\sigma)$. We will often refer to a
labeling as a coloring, but there is no assumption that adjacent vertices get
different colors. Of course the goal is to minimize the number of colors used.
Consequently  the distinguishing number of a graph $G$ is defined  by

\begin{equation*}
D(G) = min\{r \vert ~ G ~\textsl{has a labeling that is $r$-distinguishing}\}.
\end{equation*} 

This number has defined by Albertson and Collins \cite{Albert}. Similar to this definition, Kalinowski and Pil\'sniak \cite{Kali1} have defined the distinguishing index $D'(G)$ of $G$ which is  the least integer $d$
such that $G$ has an edge colouring   with $d$ colours that is preserved only by a trivial
automorphism. If a graph has no nontrivial automorphisms, its distinguishing number is  $1$. In other words, $D(G) = 1$ for the asymmetric graphs.
 The other extreme, $D(G) = \vert V(G) \vert$, occurs if and only if $G = K_n$. The distinguishing index of some examples of graphs was exhibited in \cite{Kali1}. For 
 instance, $D(P_n) = D'(P_n)=2$ for every $n\geqslant 3$, and 
 $D(C_n) = D'(C_n)=3$ for $n =3,4,5$,  $D(C_n) = D'(C_n)=2$ for $n \geqslant 6$. It is easy to see that the value $|D(G)-D'(G)|$ can be large. For example $D'(K_{p,p})=2$ and $D(K_{p,p})=p+1$, for $p\geq 4$. 
 The Cartesian product of graphs $G$ and $H$ is a graph denoted $G\Box H$ whose
 vertex set is $V (G) \times  V (H)$. Two vertices $(g, h)$ and $(g', h')$ are adjacent if
 either $g = g'$ and $hh'\in  E(H)$, or $gg' \in  E(G)$ and $h = h'$. We denote
 $G\Box G$ by $G^2$, and we recursively define the $k$-th Cartesian power of $G$ as
 $G^k = G\Box G^{k-1}$ \cite{gorz}.
A graph $G$ is called prime if $G = G_1\Box G_2$ implies that either $G_1$ or $G_2$ is
$K_1$. 
The distinguishing number and index  of the Cartesian  powers of graphs has been
thoroughly investigated. It was first proved by Albertson \cite{Albert2005} that if $G$ is a
connected prime graph, then $D(G^k) = 2$ whenever $k \geq 4$, and if $|V (G)| \geq 5$,
then also $D(G^3) = 2$. Next, Kla\v{v}zar and Zhu \cite{Klavzar}  showed that for any
connected graph $G$ with a prime factor of order at least $3$,  $D(G^k) = 2$ for $k \geq 3$. 
Michael and  Garth  in \cite{fish} have determined the distinguishing number of
the Cartesian product of complete graphs.  
Pil\'sniak studied the Nordhaus-Gaddum bounds for the distinguishing index in \cite{nord}. Also the distinguishing number of the hypercube has been investigated 
in \cite{bogs}. Similar to definition of $D(G)$ and $D'(G)$, authors in \cite{Kali2} introduced the total distinguishing number of a graph $G$, $D''(G)$  as the least number
$d$ such that $G$ has a total colouring (not necessarily proper) with $d$ colours that is only preserved
by the trivial automorphism. They proved that $D''(G)\leq \lceil\sqrt{\Delta(G)}\rceil$.  

\medskip

 In this paper, we continue the study of
 two parameters $D(G), D'(G)$ and proceed as follows.

\medskip
 In the next section, we consider two specific graphs, friendship graphs and book graphs and compute their  distinguishing number and  index. 
 Also we study the distinguishing number and the distinguishing index of  corona product of two graphs in Section 3.

\section{The distinguishing number and index of some graphs}
 In this section, we consider friendship graphs and book graphs and compute their  distinguishing number and their distinguishing  index. We begin with friendship graph. 
The friendship graph $F_n$ $(n\geqslant 2)$ can be constructed by joining $n$ copies of the cycle graph $C_3$ with a common vertex. First we state the following lemma:

\begin{figure}[ht]
	\hspace{1cm}
	\begin{minipage}{6.3cm}
		\includegraphics[width=\textwidth]{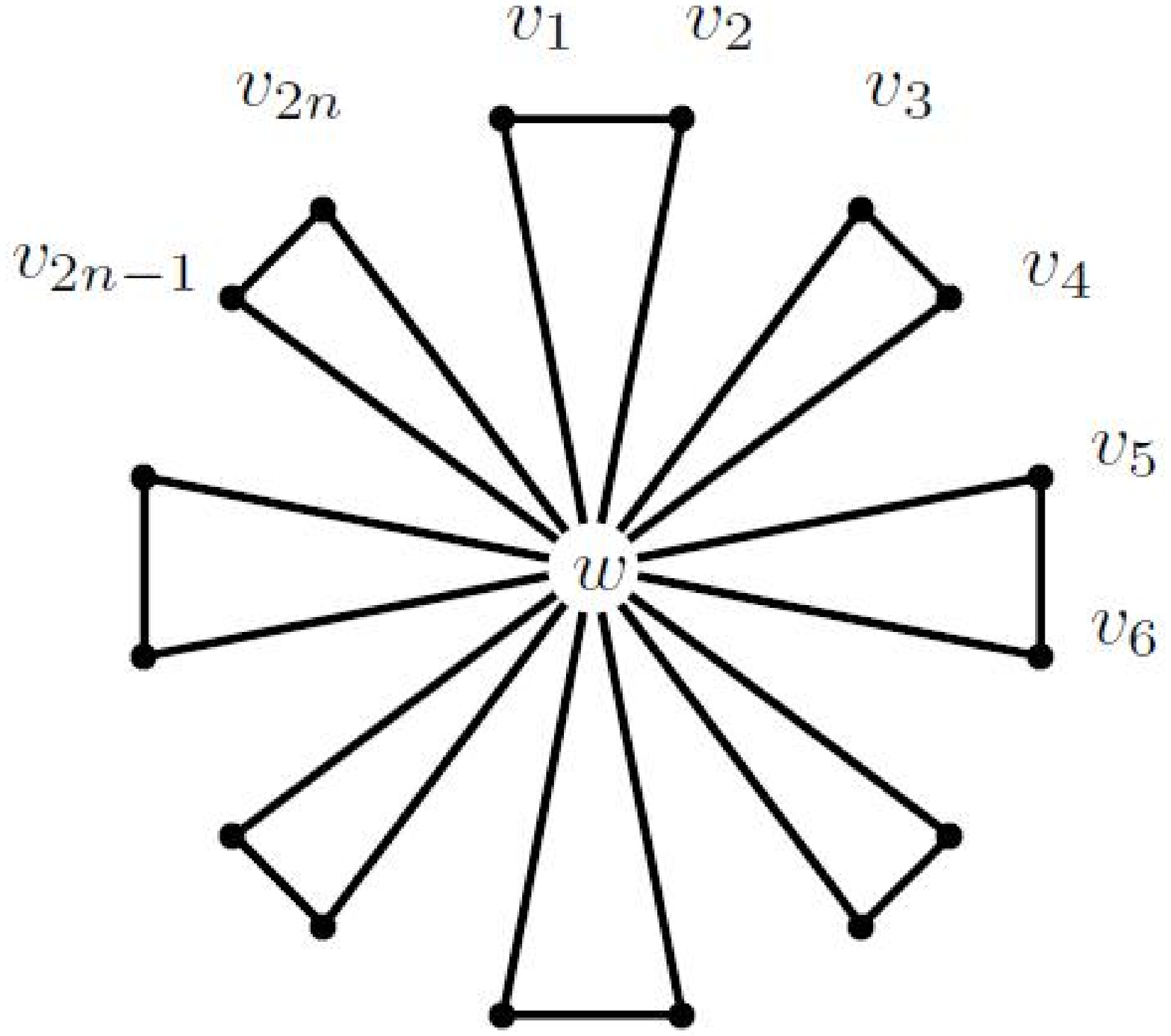}
	\end{minipage}
	\begin{minipage}{6.1cm}
		\includegraphics[width=\textwidth]{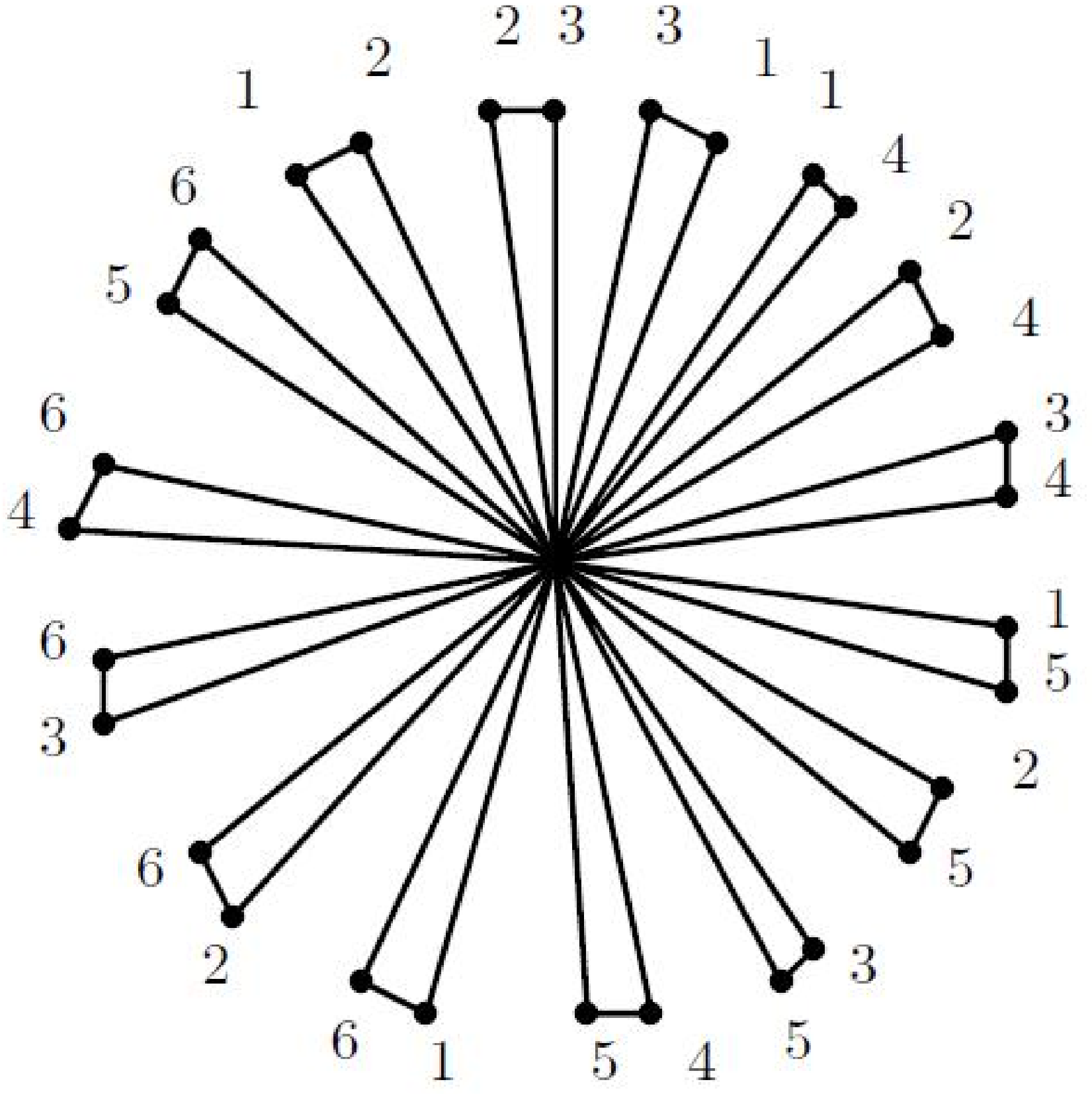}
	\end{minipage}
	\caption{\label{friend} Friendship graph $F_n$ and the vertex labeling of $F_{15}$, respectively.}
\end{figure}

\begin{lemma}
	The order of automorphism group of $F_n$   $(n\geqslant 2)$ is $|Aut(F_n)|=n!2^n$.
\end{lemma}
\proof
Since $w$ is the only vertex of $F_n$ which is not of degree $2$ (Figure \ref{friend}), so $w$ is fixed by all elements of the automorphism group. We can get  the automorphism group of $F_n$ by interchanging  the base of triangles together and rotating  about the center $w$. Therefore $|Aut(F_n)|=n!2^n$.\qed

\begin{theorem}
The distinguishing number of the friendship graph $F_n$  $(n\geq 2)$ is  $$D(F_n)= \lceil \dfrac{1+\sqrt{8n+1}}{2}\rceil .$$  
\end{theorem}
\proof
First we shall find a lower bound for $D(F_n)$ and then we present a distinguishing vertex labeling with this number of labels.
Let $\{x_i,y_i\}$, $1\leqslant i\leqslant n$ be the set of two labels that  has been assigned to the two vertices of the base of $i$-th triangle, and $L=\big\{1,\{x_i,y_i\}~\vert ~ 1\leqslant i\leqslant n,~ x_i,y_i \in \mathbb{N}\\big\}$ is the labeling of $F_n$ such that the label of the central vertex $w$ is $1$ and the label of the two vertices on the base of $i$-th triangle is $\{x_i,y_i\}$. If $L$ is a distinguishing labeling for $F_n$, then it satisfies the following properties:

\begin{itemize}
\item[(i)] For every $i\in \{1,\ldots n\}$, $x_i\neq y_i$. Because for every $1\leq i\leq n$, the map $f_i:V(F_n)\rightarrow V(F_n)$ which  maps $v_{2i-1}$ and $v_{2i}$ to each other and fixes the rest of vertices of $F_n$, is an automorphism of $F_n$.

\item[(ii)] For every $i,j\in \{1,\ldots n\}$ where $i\neq j$, $\{x_i,y_i\}\neq \{x_j,y_j\}$.  Because for every $i,j\in \{1,\ldots n\}$ where $i\neq j$, the map $f_{i,j}:V(F_n)\rightarrow V(F_n)$ which  maps $v_{2i-1}$ and $v_{2j}$ to each other and $v_{2i}$ and $v_{2j-1}$ to each other and fixes the rest of vertices of $F_n$, is an automorphism of $F_n$. Also the map $g_{i,j}:V(F_n)\rightarrow V(F_n)$ which  maps $v_{2i-1}$ and $v_{2j-1}$ to each other and $v_{2i}$ and $v_{2j}$ to each other and fixes the rest of vertices of $F_n$, is an automorphism of $F_n$. 

\end{itemize}

So it can be obtained that with labels $\{1,\ldots , s\}$ we can make at most ${s \choose 2}$ numbers of   the pairs $(x,y)$  such that they satisfy (i) and (ii). Hence   $D(F_n) \geqslant min\{s: {s \choose 2}\geqslant n \}$. By a simple computation we get $D(F_n) \geqslant \lceil \dfrac{1+\sqrt{8n+1}}{2}\rceil $. 
Now we  define a distinguishing vertex labeling on $F_n  $ with $\lceil \dfrac{1+\sqrt{8n+1}}{2}\rceil$  labels. Consider the friendship graph in Figure \ref{friend}. The function that maps $ v_1 $ to $ v_2 $ and $ v_2 $ to $ v_1 $ and fixes the rest of vertices, is a non-trivial automorphism. Thus the labels $ v_1$ and $v_2 $ should be different.  We assign the vertex $ v_1 $  the label $1$ and the vertex $ v_2$ the label $2$. Similarly, the function that maps $ v_3 $ to $ v_4 $ and $ v_4 $ to $ v_3 $ and fixes the rest, is a non-trivial automorphism. Thus the labels $ v_3 , v_4 $ should be distinct. Let assign the vertex $ v_3 $ the label $2$  and the vertex  $ v_4$ the label $3$. We continue this method to label all vertices of friendship graph (see the label of $F_{15}$ in Figure \ref{friend}). Note that the label of vertex $w$ is $1$. 
Hence this  method gives a distinguishing vertex labeling with the minimum number of labels.
By the above process, observe that the distinguishing number of  $F_n$, $D(F_n)$, is the $n$-th term of the sequence $\{D(F_i) \}$ which defines as follows:  
\begin{equation*}
\{D(F_i) \}_{i\geqslant 1}=\{-,3,3,4,4,4,\underbrace{5,\cdots ,5}_{4-times},\underbrace{6,\cdots ,6}_{5-times},\underbrace{7,\cdots ,7}_{6-times},\cdots ,\underbrace{m,\cdots ,m}_{(m-1)-times}, \cdots \}.
\end{equation*}
In fact,
\begin{equation*}
D(F_n)  =min \{k: \sum_{i=2}^k (i-1)\geqslant n\}.
\end{equation*}
By an easy computation, we see that 
$$min \{k: \sum_{i=2}^k (i-1)\geqslant n\}=\lceil \dfrac{1+\sqrt{8n+1}}{2}\rceil.$$
Therefore we have the result.\qed

\begin{figure}
	\begin{center}
		\includegraphics[width=0.5\textwidth]{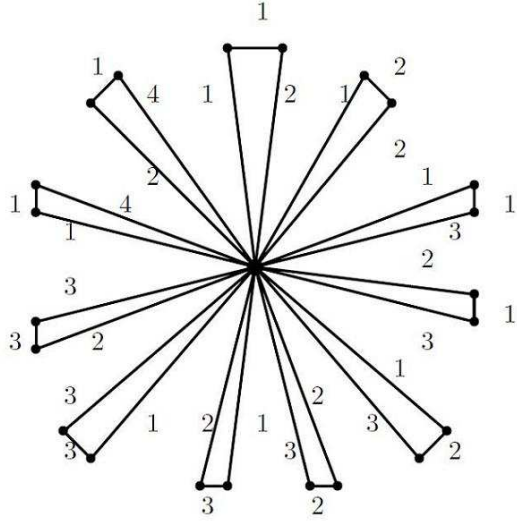}
		\caption{The edge labeling of $F_{11}$.}\label{edge-label friend}
	\end{center}
\end{figure}

Now we shall compute the distinguishing index of the friendship graph, i.e., $D'(F_n)$. 
\begin{theorem}
Let $a_n=1+27n+3\sqrt{81n^2+6n}$. 
	For every $n\geq 2$, $$D'(F_n)=\lceil\frac{1}{3} (a_n)^{\frac{1}{3}}+\frac{1}{3(a_n)^{\frac{1}{3}}}+\frac{1}{3}\rceil.$$
\end{theorem} 
\proof
First we show that  $D'(F_n)\geqslant min\{k: k^3-k^2\geq 2n\}$ and next we present a distinguishing edge labeling such that it obtains this bound.

Let $\{x_i,y_i, z_i\}$, $1\leqslant i\leqslant n$ be the label of three sides of a triangle in the friendship graph  such that $ z_i $ is the label which it is assigned to the base and $ x_i,y_i $ are the labels of two sides,  and let $L'=\{(x_i,y_i,z_i)~\vert ~ 1\leqslant i\leqslant n,~ x_i,y_i \in \mathbb{N}\}$ be the labeling of $F_n$. If $L'$ is a distinguishing labeling for $F_n$, then it satisfies the following properties:

\begin{itemize}
\item[(i)] For all $j=1,\ldots, n$, $(x_j,y_j,z_j) \neq (y_j, x_j, z_j)$. Because for every $i\in \{1,\ldots n\}$, the map $f_i:V(F_n)\rightarrow V(F_n)$ which maps $v_{2i-1}$ and $v_{2i}$ to each other and fixes the rest of the vertices of $F_n$, is an automorphism of $F_n$.

\item[(ii)] For every $j\neq i$,   $(x_i,y_i,z_i) \neq (x_j, y_j, z_j)$ and $(y_i,x_i,z_i) \neq (x_j, y_j, z_j)$. Because for every $i,j\in \{1,\ldots n\}$ where $i\neq j$, the map $f_{i,j}:V(F_n)\rightarrow V(F_n)$ which  maps $v_{2i-1}$ and $v_{2j}$ to each other and $v_{2i}$ and $v_{2j-1}$ to each other and fixes the rest of vertices of $F_n$, is an automorphism of $F_n$. Also the map $g_{i,j}:V(F_n)\rightarrow V(F_n)$ such that it maps $v_{2i-1}$ and $v_{2j-1}$ to each other and $v_{2i}$ and $v_{2j}$ to each other and fixes the rest of vertices of $F_n$, is an automorphism of $F_n$.
\end{itemize}

So it can be obtained that with labels $\{1,\ldots , s\}$ we can make at most ${s \choose 2}s$  numbers of   the $3$-ary's $(x,y,z)$  such that they satisfy (i) and (ii) (there are ${s \choose 2}$ choices for $x_i$ and $y_i$ and $s$ choices for $z_i$). Hence  $D'(F_n) \geqslant min\{s:  {s \choose 2}s\geq n \}$ and it can be calculated that $D'(F_n) \geqslant \lceil\frac{1}{3} (a_n)^{\frac{1}{3}}+\frac{1}{3(a_n)^{\frac{1}{3}}}+\frac{1}{3}\rceil$.

Now we  define a distingushing edge labeling on $F_n  $ with $ min\{k: k^3-k^2\geq 2n\}$ labels. Similar to the vertex labeling of $F_n$, in the edge labeling of $F_n$,  the labels of two sides of every triangle should be distinct, otherwise, we have  a non-trivial automorphism, which preserves the labeling.  We assign the first triangle,  the $3$-ary  $(1,2,1)$ and the second  the $3$-ary  $(1,2,2)$.  Now we assign the third  triangle,  the $3$-ary  $(1,3,1)$  and the forth  triangle  the $3$-ary  $(2,3,1)$. Continuing this method we can obtain a distinguishing labeling for the graph (see the Figure \ref{edge-label friend}  for the labeling of $F_{11}$).   
It is easy to see that the distinguishing index of  $F_n$, $D'(F_n)$, is the $n$-th term of the sequence $\{D'(F_i) \}$ which defines as follows:  
\begin{equation*}
\{D'(F_i) \}_{i\geqslant 1}=\{-,2,\underbrace{3,\ldots , 3}_{7-times},\underbrace{4, \ldots ,4}_{15-times},\underbrace{5,\ldots ,5}_{26-times},\underbrace{6,\ldots ,6}_{40-times},\ldots ,\underbrace{m,\ldots ,m}_{2(m-1)+3{m-1\choose 2}-times},\ldots \}.
\end{equation*}
In fact,
\begin{equation*}
 D'(F_n)  =min \{k: \sum_{i=2}^k \big(2(i-1)+3{i-1\choose 2} \big)\geqslant n\}.
\end{equation*}
By an easy computation, we see that 
\begin{align*}
&~~~~~~~min \{k: \sum_{i=2}^k \big(2(i-1)+3{i-1\choose 2} \big)\geqslant n\}=min\{k: k^3-k^2\geq 2n\}=\\
& \lceil \frac{1}{3}(1+27n+3\sqrt{81n^2+6n})^{1/3}+\dfrac{1}{3(1+27n+3\sqrt{81n^2+6n})^{1/3}} +\frac{1}{3}\rceil.
\end{align*}
So, our method for edge labeling of $F_n$ which as shown in Figure \ref{edge-label friend} lead to use the minimum number of labels. Therefore we have the result.\qed

\begin{figure}
	\begin{center}
		\hspace{.7cm}
		\includegraphics[width=0.5\textwidth]{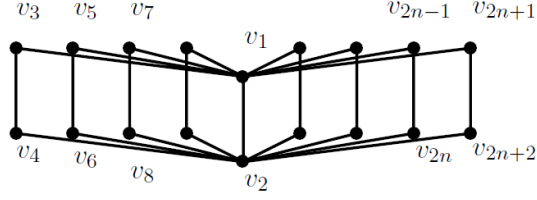}
				\caption{Book graph $B_n$.}\label{book} 
	\end{center}
\end{figure}

 The $n$-book graph $(n\geqslant 2)$ (Figure \ref{book}) is defined as the Cartesian product $K_{1,n}\square P_2$. We call every $C_4$ in the book graph $B_n$, a page of $B_n$. All pages in $B_n$ have a common side $v_1v_2$.  If we change the labels  of vertices of    parallel side  of $v_1v_2$ (for example the labels of $v_3$ and $v_4$ in Figure \ref{book}), then we call this new page as the inverse of the first page. We shall compute the distinguishing number and index of $B_n$. The following result gives the order  of automorphism group  $B_n$.

\begin{theorem}
	For every $n\geq 2$, $|Aut(B_n)|=2n!$.
\end{theorem}
\proof 
All vertices of $B_n$, except two vertices  $v_1$ and $v_2$ have degree $2$ (Figure \ref{book}). So the two vertices $v_1$ and $v_2$ are mapped into each other under the elements of automorphism group. In fact each automorphism maps the pages to each other. Note that as soon as the first page is mapped to the inverse of a page, the rest of pages are mapped to inverse of themselves or to inverse of another page. Therefore $\vert Aut(B_n)\vert = 2n!$.\qed

\begin{theorem}
	The distinguishing number of $B_n$ $(n\geq 2)$ is $D(B_n)= \lceil \sqrt{n}\rceil $.
\end{theorem}
\proof First we show that  $D(B_n)\geqslant \lceil \sqrt{n}\rceil$ and next we present a distinguishing edge labeling such that it obtains this bound.
Let $x_i, y_i$ be the labels of upper and lower vertices of $i$-th page respectively (except for $v_1,v_2$ in Figure \ref{book})  and let $L=\{(1,2), (x_i,y_i)~\vert ~ 1\leqslant i\leqslant n,~ x_i,y_i \in \mathbb{N}\}$ be the labeling of $B_n$ such that $(1,2)$ is the two labels that has been assigned to the vertices $v_1$ and $v_2$. If $L$ is a distinguishing labeling for $B_n$, then it  has the following property:
\begin{itemize}
\item[(i)] For all $i,j\in\{1,\ldots , n\}$ where $j\neq i$,~~$(x_i,y_i)\neq (x_j,y_j)$,
\end{itemize} 
because for every $i,j\in \{1,\ldots n\}$, the map $f_i:V(B_n)\rightarrow V(B_n)$ which  maps $v_{2i-1}$ and $v_{2j-1}$ to each other and maps $v_{2i}$ and $v_{2j}$ to each other and fixes the rest of vertices of $B_n$, is an automorphism of $B_n$.

So it can be obtained that with labels $\{1,\ldots , s\}$ we can make at most $2{s \choose 2}+s$ numbers of  the pairs $(x,y)$  such that they satisfy (i) (there are    $2{s \choose 2}$ choices for  pairs $(x,y)$ such that $x\neq y$ and $(x,y)$ satisfies (i) and  since the pairs $(i,i)$ for  every $i\in \{1,\ldots , s\}$ are not counted in these $2{s \choose 2}$ choices, so we add $s$ to $2{s \choose 2}$). Hence  $D(B_n) \geqslant min\{s:  2{s \choose 2}+s\geq n \}$ and it can be calculated that $D(B_n) \geqslant  \lceil \sqrt{n}\rceil$.

\begin{figure}
	\begin{center}
		\includegraphics[width=0.4\textwidth]{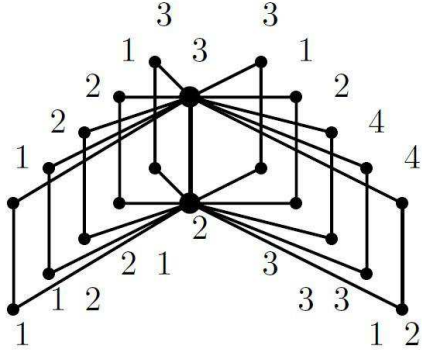}
		\caption{The vertex labeling of $B_{10}$.}\label{v-l book} 
	\end{center}
\end{figure}

To present  the vertex labeling of $B_n$, note that the labels of $v_1$ and $v_2$ can be the same or distinct. Since we would like to use least number of labels, observe that for this purpose,  two labels that have been given to $v_1$ and $v_2$ should be different, because if the label of vertices   $v_1$ and $v_2$ is the same, then we should also check $x_i\neq y_i$ for each $i$ $(1\leqslant i\leqslant n)$. We assign the first page of $B_n$, the pair  $(1,1)$ and the second,  the pair  $(2,1)$.  Now we assign the third  page of $B_n$,  the pair  $(2,2)$  and the forth  page,  the pair $(1,2)$. Now we use the new label $3$ for the labeling next five pages. We assign these five pages the labels $(3,1)$, $(3,2)$, $(3,3)$, $(1,3)$, and $(2,3)$, respectively. Note that if in the process of labeling,  the same labels have been given to $v_1$ and $v_2$, then commute them with two labels of the  next page.  Our method for labeling the vertices of book graph have been shown for $B_{10}$ in Figure \ref{v-l book}.
By this  process, observe that the distinguishing number of $B_n$, $D(B_n)$, is the $n$-th term of the sequence $\{D(B_i)\}$ which defines as follows: 
\begin{equation*}
\{D(B_i)\}_{i\geqslant 1}=\{-,2,2,2,3,3,3,3,3,\underbrace{4,\cdots ,4}_{7-times},\underbrace{5,\cdots ,5}_{9-times},\cdots ,\underbrace{m,\cdots , m}_{(2m-1)-times}, \cdots \}.
\end{equation*}
In fact, 
\begin{equation*}
D(B_n)=min \{k: \sum_{i=1}^k (2i-1)\geqslant n\}.
\end{equation*}
By an easy computation, we see that

\begin{equation*}
min \{k: \sum_{i=1}^k (2i-1)\geqslant n\}= \lceil \sqrt{n}\rceil.
 \end{equation*}

Therefore we have the result.\qed

\begin{remark}
	  The distinguishing index of  Cartesian product of star $K_{1,n}$ with path $P_m$ for $m \geqslant≥ 2$ and $n \geqslant 2$ is $D'(K_{1,n}\square P_m) = \lceil \sqrt[2m-1]{n} \rceil$, unless $m=2$ and $n=r^3$ for some integer $r$. In the latter case $D'(K_{1,n}\square P_2) = \sqrt[3]{n}+1$. (\cite{gorz}).  Since $B_n=K_{1,n}\Box P_2$,  using this equality  we obtain the distinguishing index of book graph $B_n$.
\end{remark}

\section{Distinguishing number (index) of corona of two graphs}
In this section, we shall study the distinguishing number and the distinguishing index  of corona product of two graphs.  
The corona product $G\circ H$ of two graphs $G$ and $H$ is defined as the graph obtained by taking one copy of $G$ and $\vert V(G)\vert $ copies of $H$ and joining the $i$-th vertex of $G$ to every vertex in the $i$-th copy of $H$.

\begin{theorem}
	For every $n\geq 4$, $D(P_n\circ K_1)=(D'(P_n\circ K_1))=2$.
	\end{theorem}
\proof
Since $ \vert Aut(P_n\circ K_1 ) \vert =2$, so $ P_n\circ K_1 $ has only one non-trivial automorphism. Therefore $ D(P_n\circ K_1)=D'(P_n\circ K_1)=2$.\qed

Before presenting the main result for   $ D(G\circ H) $, we explain the relationship between the automorphism group of the graph $ G\circ H $ with the automorphism groups of two connected  graphs $ G$ and $H$ such that $G\neq K_1$.
 Note that  there is no vertex in the copies of $H$ which has the same degree as a vertex in $G$. Because if  there exists a vertex $w$ in one of the copies of $H$  and a vertex $v$ in $G$ such that $deg_{G\circ H} v=deg_{G\circ H} w$, then
    ${deg}_G(v)+|V(H)|={deg}_H(w)+1$. So we have $deg_H(w)+1 > |V(H)|$, which is a contradiction.
   By this note, we state and prove  the following theorem: 
   
\begin{theorem}\label{thh1}
	For every two connected graphs $G$ and $H$ such that $G\neq K_1$, we have
$\vert Aut (G\circ H)\vert = \vert Aut (G)\vert  \vert Aut(H)\vert$.
\end{theorem}
\proof
 Let the vertex set of  $G$ be $\{v_1,\ldots , v_{\vert V(G)\vert}\}$ and the vertex set of $ i$-th copy of $H$, $H^{(i)}$, be $\{w^{(i)}_1,\ldots , w^{(i)}_{\vert V(H)\vert}\}$. 
Since there is no vertex in copies of $H$ which has the same degree as a vertex in $G$,  for every $f \in Aut( G\circ H )$, we have 
$f\vert _H\in Aut(H)$ and $f\vert _G\in Aut(G).$
In addition, for $i,j\in \{1,\ldots , \vert V(G)\vert\}$ we have

\begin{equation*}
f(v_i)=v_j \Longleftrightarrow f(H^{(i)})=H^{(j)}.
\end{equation*}

Conversely, let $ \varphi \in Aut(G) $ and  $ \phi \in Aut(H) $ such that $\varphi (v_i)=v_{j_i}$, where $i,j_i\in \{1,\ldots , \vert V(G)\vert\}$. Now we define the following automorphism $ h $ of $G\circ H$:
\begin{equation*}
\left\{
\begin{array}{ll}
h: G\circ H\rightarrow G\circ H & \\
~~~v_i\mapsto \varphi(v_i)=v_{j_i} &~ i,j_i\in \{1,\ldots , \vert V(G)\vert\},\\
~~~w_{k}^{(i)}\mapsto (\phi(w_k))^{(j_i)}&~ k\in \{1,\ldots , \vert V(H)\vert\}.
\end{array}\right.
\end{equation*}

Therefore $\vert Aut (G\circ H)\vert = \vert Aut (G)\vert  \vert Aut (H)\vert$.\qed

   \medskip

 \begin{figure}
 	\begin{center}
 		\includegraphics[width=0.9\textwidth]{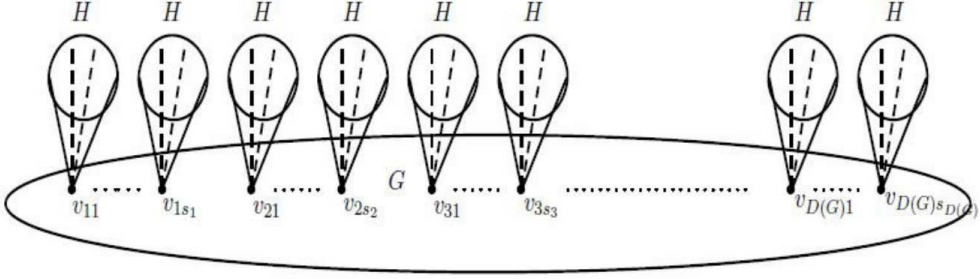}
 		\caption{The partition of the vertices of $G\circ H$ by its labels.}\label{GoH}
 	\end{center}
 \end{figure}
 
By  the elements of the automorphism group of $G\circ H$ (Theorem \ref{thh1}), we have  the following result.  
\begin{theorem}\label{thh2} Let $G$ and $H$ be two connected graphs and $G\neq K_1$.
\begin{enumerate}
\item[(i)] If $ D(G)=1 $, then $D(G\circ H) =D(H)$.\label{1}
\item[(ii)] $D(G\circ H) =1$ if and only if $D(G)=D(H)=1$.\label{3}
\end{enumerate}
\end{theorem}

\begin{theorem}\label{specific}
	Let $G$ and $H$ be two connected graphs such that $G\neq K_1$. If  $D(G)\leqslant D(H)$, then $D(H)= D(G\circ H)$
\end{theorem}
\proof
If $D(G)=1$, then we have the result by Theorem \ref{thh2}, so we suppose that $D(G)\neq 1$.
 If we label $G\circ H$ with less than $D(H)$ labels in a distinguishing way, then we can find a non-identity automorphism of $H$ such as $f$, such that it preserves the labeling of $H$. Expanding $f$ to $G\circ H$ such that $f$ acts as the identity function on $G$, we obtain a non-identity automorphism of $G\circ H$  preserving the labeling of $G\circ H$, which is contradiction. So we have  $D(H)\leqslant D(G\circ H)$. 
Now we show the inequality  
$D(H)\geqslant D(G\circ H)$.  By the  definition of distinguishing vertex labeling, the vertex set $V(G)$ is  partitioned  to  at most $D(H)$-classes (because $D(G)\leqslant D(H)$), say, $[1],[2],...,[D(H)]$. The vertices of the class $[i]$ denoted  by $v_{i1},\ldots , v_{is_i}$  in  Figure \ref{GoH} where $s_i$ is the size of the $[i]$-class and $ i=1,\ldots , D(H)$. We label the vertices  in the class $[i]$ and $s_i$-copies of $H$ to obtain a distinguishing vertex labeling of   $G\circ H$ as follows: we label  all vertices in the class $[i]$  with the label $i$, and vertices in the $s_i$ copies of $H$ with $D(H)$ labels in a distinguishing way where $ i=1,\ldots , D(H)$. By Theorem \ref{thh1} this labeling is a distinguishing labeling of $G\circ H$, and so $D(H)\geqslant D(G\circ H)$.\qed

\begin{theorem}  \label{special}
	Let $G$ and $H$ be two connected graphs such that  $ D(G)> D(H)$.  Then 
\begin{equation*}
D(H)\leqslant D(G\circ H) \leqslant D(H)+\lfloor \dfrac{-(1+D(H))+\sqrt{(D(H)-1)^2 +4D(G)}}{2}  \rfloor.
\end{equation*}
\end{theorem}
\proof
The proof of the left inequality is the same as the first part of the proof of Theorem \ref{specific}. For the second inequality, by the  definition of distinguishing vertex labeling, the vertex set $V(G)$ will be partitioned  to  $D(G)$-classes, say, $[1],[2],...,[D(G)]$. The vertices of the class $[i]$ denoted  by $v_{i1},\ldots , v_{is_i}$ ($ i=1,\ldots , D(G)$) in  Figure \ref{GoH}. We label the vertices  in the class $[i]$ and $s_i$-copies of $H$ to obtain a distinguishing vertex labeling of   $G\circ H$ as follows: 

{\bf Step 1)} Labeling  the vertices in the  classes $[i]$  and the vertices of $s_i$-copies of  $H$ for $1\leq i\leq D(H)$: 

We label  all vertices in the class $[i]$  with the label $i$, and vertices in the $s_i$ copies of $H$ with $D(H)$ labels in a distinguishing way. 

\medskip
{\bf Step 2)}  Labeling the vertices in the  classes $[i]$ and the vertices of $s_i$-copies of  $H$, for  $D(H)+1\leq i\leq 2D(H)+2$: 

Now for the labeling of the vertices of classes $[i]$, $(D(H)+1\leq i\leq 2D(H)+1)$ we use the label $i-D(H)$, and for $s_i$ copies of $H$, we add the number one to the label of each vertex of $H$ in the prior step, i.e., if the label of a vertex of $H$ is $l$ $(1\leqslant l\leqslant D(H))$ in distinguishing labeling of $H$ with $D(H)$ labels, then we replace it by $l+1$. For the class $[2D(H)+2]$ we use the label $D(H)+1$ for the vertices in this class and we label $s_{2D(H)+2}$ copies of $H$ with $D(H)$ labels in a distinguishing way.
 \medskip

{\bf Step 3)}  Labeling the vertices in the  classes $[i]$ and the vertices of $s_i$-copies of  $H$, for  $2D(H)+3\leq i\leq 3D(H)+6$: 

Now for the labeling of the vertices of classes $[i]$, $(2D(H)+3\leq i\leq 3D(H)+4)$ we use the label $i-(2D(H)+2)$, and for $s_i$ copies of $H$, we add the number two to the label of each vertex of $H$ in the firs step, i.e., if the label of a vertex of $H$ is $l$ $(1\leqslant l\leqslant D(H))$ in distinguishing labeling of $H$ with $D(H)$ labels, then we replace it by $l+2$. For the class $[3D(H)+5]$ we use the label $D(H)+2$ for the vertices in this class and label $s_{3D(H)+5}$ copies of $H$ with $D(H)$ labels in a distinguishing way. For the class $[3D(H)+6]$ we use the label $D(H)+2$ for the vertices in this class and we label $s_{3D(H)+6}$ copies of $H$ with $D(H)+1$ labels in a distinguishing way (i.e.,  if the label of a vertex of $H$ is $l$ $(1\leqslant l\leqslant D(H))$ in distinguishing labeling of $H$ with $D(H)$ labels, then we replace it by $l+1$).

Continuing this method and by Theorem \ref{thh1} it can be observed that  this method makes a distinguishing labeling with $D(H)+min\{k: \left( \sum_{i=0}^{k} (D(H)+2i)\right) \geqslant D(G)\}$ labels. By an easy computation we get
\begin{equation*}
min\{k: \left( \sum_{i=0}^{k} (D(H)+2i)\right) \geqslant D(G)\} =\lfloor \dfrac{-(1+D(H))+\sqrt{(D(H)-1)^2 +4D(G)}}{2}  \rfloor.
\end{equation*}

So  we have the result.\qed

\begin{theorem}
Let $H$ be a conneced graph, then $D(H)\leqslant D(K_1\circ H)\leqslant D(H)+1$.
\end{theorem}
\proof
First we prove $D(H)\leqslant D(K_1\circ H)$.  Suppose to the contrary that  $D(H)> D(K_1\circ H)$, so if we label $K_1 \circ H$ in a distinguishing way with $D(K_1\circ H)$ labels and transfer this labeling to $H$, then there exists a non-identity automorphism of $H$ such as  $f$, that it preserves the labeling. Expanding $f$ to $K_1\circ H$ such that $f$ acts as the identity on $K_1$, we have a non-identity automorphism of $K_1\circ H$ that it preserves the labeling, which is a contradiction. Now we shall show that 
$D(K_1\circ H)\leqslant D(H)+1$.  For this purpose, we define a distinguishing labeling of $K_1\circ H$ with $D(H)+1$ labels. First we label $H$  with $D(H)$ labels in a distinguishing way and next assign a new label to the only vertex of $K_1$. This labeling is a distinguishing labeling for $K_1\circ H$, because if $f$ is an automorphism of $K_1\circ H$ preserving the labeling, then $f (K_1) =K_1$ and $f\vert _H \in Aut(H)$. Since we labeled $H$ in a distinguishing way, $f\vert_H$ is the identity automorphism. Therefore $f$ is the identity  automorphism on $K_1\circ H$. Therefore, the result follows.\qed

Here we study the distinguishing index of corona of two graphs. First we compute the distinguishing index of some special cases and   exclude them subsequently. The special cases are as follows:
\begin{equation*}
D'(K_1\circ K_1)=1, ~~D'(K_1\circ K_2)=3,~~ D'(K_2\circ K_1)=2,~~ D'(K_2\circ K_2)=2.
\end{equation*}

\begin{theorem}
Let $G$ and $H$  be two connected graphs such that $G\neq K_1$ and $D'(H)\geqslant 2$, then $D'(G\circ H) \leqslant max\{D'(G), \lceil \sqrt{D'(H)} \rceil\}$.
\end{theorem}
\proof
We define a distinguishing edge labeling for $G\circ H$ with $max\{D'(G), \lceil \sqrt{D'(H)} \rceil\}$ labels. First we label $G$  with the labels $\{1,\ldots , D'(G)\}$ in a distinguishing way. Now we present a labeling for a copy of $H$ and all middle edges that are incident to this copy of $H$ and $G$, and next we transfer this labeling to all copies of $H$ and their middle edges. For this we partition the edge set of $H$ with respect to a distinguishing edge labeling of $H$ with the label set $\{1, \ldots , D'(H)\}$. So we have $D'(H)$ classes of edges such that $[i]$-class $(1\leqslant i \leqslant D'(H))$ contains  all the edges of $H$ which they have the label $i$ in the distinguishing edge labeling of $H$. It is clear that there are vertices of $H$ that  are incident to the edges in different classes, such as $[i]$ and $[j]$ with $i\geqslant j$. In this case the middle edges incident to such vertex are considered as the middle edges of the $[i]$-class. The new labeling of $H$ and all its middle edges are as follows:

\medskip
{\bf Step 1)} We label all edges in class $[1]$ with the label $1$. Next we label all its middle edges  that are incident to a vertex in $[1]$-class, with the label $1$. We label all edges in class $[2]$ with the label $1$. Next we label all its middle edges  that are incident to a vertex in $[2]$-class, with the label $2$.

 \medskip
{\bf Step 2)} We label all edges in class $[3]$ with the label $2$. Next we label all its middle edges  that are incident to a vertex in $[3]$-class with the label $1$.  We label all edges in class $[4]$ with the label $2$. Next we label all its middle edges  that are incident to a vertex in $[4]$-class with the label $2$.

\medskip
{\bf Step 3)}   We label all edges in class $[5]$ with the label $1$. Next we label all its middle edges  that are incident to a vertex in $[5]$-class with the label $3$.  We label all edges in class $[6]$ with the label $2$. Next we label all its middle edges  that are incident to a vertex in $[6]$-class with the label $3$.  We label all edges in class $[7]$ with the label $3$. Next we label all its middle edges   that are incident to a vertex in $[7]$-class with the label $3$. 

 \medskip
{\bf Step 4)} We label all edges in class $[8]$ with the label $3$. Next we label all its middle edges  that are incident to a vertex in $[8]$-class with the label $1$.  We label all edges in class $[9]$ with the label $3$. Next we label all its middle edges  that are incident to a vertex in $[9]$-class with the label $2$.

Continuing this method, in the next step we label all edges in class $[10]$ with the label $4$ and next  we label all its middle edges  that are incident to a vertex in $[10]$-class with the label $1$, we  obtain a labeling for $G\circ H$ that is distinguishing. Because if $f$ is an automorphism of $G\circ H$ preserving the labeling, then the restriction of $f$ to $G$ is the identity automorphism of $G$. On the other hand, for each non-identity automorphism of $H$, there exists an edge in a class that is mapped to an edge in another class. So by considering our labeling of $H$ and all its middle edges we obtain that the restriction of $f$ to $H$ is the identity automorphism of $H$. Therefore $f$ is the identity automorphism of $G\circ H$. Since we used $min\{k: \sum_{i=1}^k (2i-1)\geqslant D'(H)\}$ labels for the labeling of copies of $H$ (and since  this number is equal with  $\lceil  \sqrt{D'(H)}\rceil$) and used $D'(G)$ labels for $G$, so we have the result. \qed

\begin{theorem}
Let $G$ and $H$ be two connected graphs of  orders  $n \geqslant3$ and $m \geqslant3$, respectively. If $D'(G)=D'(H)=1$ then $D'(G\circ H)=1$.
\end{theorem}
\proof Since  the orders of $G$ and $H$  are greater than two and $D'(G)=D'(H)=1$, so $\vert Aut (G) \vert= \vert Aut (H) \vert =1$. By Theorem \ref{thh1}, $\vert Aut (G\circ H) \vert =1$, and so $D'(G\circ H)=1$.\qed

\begin{theorem}
Let  $H$ be a connected graph of order $n \geqslant3$. Then $D'(K_1\circ H)\leqslant D'(H)+1$. 
\end{theorem}
\proof   We label the edges of $H$ with the labels $\{1,\ldots , D'(H)\}$ in a distinguishing way and next label all its middle edges with the new label $0$. If $f$ is an automorphism of $K_1\circ H$ preserving the labeling, then $f(K_1)=K_1$ and $f\vert_H$ is an automorphism of $H$. Since we labeled $H$ in a distinguishing way, so this labeling is a distinguishing labeling for $K_1\circ H$. Hence $D'(K_1\circ H)\leqslant D'(H)+1$. \qed

\begin{theorem}
Let $G$ be a connected graph such that $G\neq K_1$. Then $D'(G\circ K_2)\leqslant max\{D'(G),2\}$.
\end{theorem}
\proof   If we label $G$ with $D'(G)$ labels in a distinguishing way and label all copies of $K_2$ with the label $1$ and next assign the two middle edges  of each copy of $K_2$, the labels $1$ and $2$, then we have a distingushing labeling of $G\neq K_1$ with $max\{D'(G),2\}$ labels. \qed

\begin{theorem}\label{circ}
Let $G$ and $H$ be two connected graphs such that $G\neq K_1$ and $H\neq K_2$.
\begin{enumerate}
\item[(i)] If $\vert Aut(G) \vert =1$, then $ D'(G\circ H)\leqslant min\{D'(H),\vert V(H)\vert \}$.
\item[(ii)] If $\vert V(G)\vert \leqslant \vert V(H)\vert +1 $ and $D'(H)=1$, then $D'(G\circ H)\leqslant 2$.
\end{enumerate}
\end{theorem}
\proof
\begin{enumerate}
	\item[(i)] 
	If $ \vert Aut(G) \vert =1$, then every element of the automorphism group of $G\circ H$ treats as the identity on $G$. If $ \vert V(H)\vert < D'(H)$ then we assign the edges  between $G$ and $H^{(i)}$, the labels $1,2,\ldots , \vert V(H)\vert $ for $1\leqslant i\leqslant \vert V(G)\vert $ and assign the remaining edges the label $1$. If $ \vert V(H)\vert \geqslant D'(H)$, then we label each copy of $H$ with $D'(H)$ labels in a distinguishing way and assign the remaining edges the label $1$. In both cases we made a distinguishing edge labeling, and so the result follows.

\item[(ii)] 
Let the vertex set of  $G$ be $\{v_1,\ldots , v_{\vert V(G)\vert}\}$ and the vertex set of $ i$-th copy of $H$ be $\{w^{(i)}_1,\ldots , w^{(i)}_{\vert V(H)\vert}\}$.  Let $ e_{ik} $ be the edge from $v_i$ to $w_k^{(i)}$.  If $ \vert Aut(G) \vert \geqslant 2$, then there exists a non-trivial automorphism $ \varphi $ of $ G\circ H $ and $ r,s\in \{1,\ldots,\vert V(G)\vert\}$, $ r\neq s $ such that $ \varphi(v_r)=v_s$. So $ e_{rk} $ is mapped to $ e_{sk'} $ under $ \varphi $ where $k,k' \in \{1,\ldots , \vert V(H)\vert \}$. Now we assign $ e_{i1},\ldots , e_{i(i-1)} $ the label $2$ for $  2\leqslant i \leqslant \vert V(G) \vert$, and assign the remaining edges the label $  1$. Clearly, this labeling is distinguishing (see Theorem \ref{thh1}), and so $D'(G\circ H)\leqslant 2$.
\qed
\end{enumerate}
\begin{corollary}
Let $G$ and $H$ be two connected graphs such that $G\neq K_1$  and $H\neq K_2$.
\begin{enumerate}
\item[(i)] If $D(G) =1$, then  $D'(G\circ H)\leqslant min\{D'(H), \vert V(H)\vert \}$.
\item[(ii)] If $\vert V(G)\vert \leqslant \vert V(H)\vert +1 $ and $D'(H)=1$, then $D'(G\circ H)\leqslant 2$.
\end{enumerate}
\end{corollary}
\proof
\begin{enumerate}
\item[(i)] 
We note that for every graph $G$,  $D(G) =1$ if and only if  $\vert Aut(G) \vert =1$. So we have the result by Theorem \ref{circ} (i). 
\item[(ii)] 
It is easy to see that for every  $G$, $D(G) \geqslant 2$ if and only if $\vert Aut(G) \vert \geqslant 2$.  So we have the result by Theorem \ref{circ} (ii).\qed
\end{enumerate}

Now, we shall present an upper bound for $D'(G\circ H)$ with $D'(H)=1$ without any condition on $\vert V(G)\vert$. For this purpose we need two following parameters:
\begin{equation*}
x'_r=\left\{
\begin{array}{ll}
 1&r=1\\
 m-1&r=2\\
 \sum_{i_{r-2}=r-1}^{m}\cdots \sum_{i_2=i_3}^{m}\sum_{i_1=i_2}^{m}(m-i_1)&r\geqslant 3,
 \end{array}\right. 
\end{equation*}
\begin{equation*}
 y'_r=\left\{
 \begin{array}{ll}
 1&r=1\\
 m&r=2\\
 \sum_{i=0}^{r-1}{r-1\choose i}x_{i+1}&r\geqslant 3.
\end{array}\right.
\end{equation*}

In fact, $x'_r$ is the number of copies of $H$ in $G\circ H$,  that their middle edges (edges between  $H$ and $G$),  have been labeled with  $r$ labels such that these   $r$ labels are used in each copy at least one time. Also $y'_r$ is the number of copies of $H$ that their middle edges, the edges between $H$ and $G$,    have been labeled with the labels $1,\ldots , r$ such that the label $r$ is used in each copy at least one time.
\begin{theorem}
Let $G$ and $H$ be the two connected graphs of orders $n$ and $m$, respectively such that  $G\neq K_1$ and $H\neq K_2$ and $D'(H)=1$. If $D(G)\geqslant 2$ then $D'(G\circ H)\leqslant min\{D'(G), min\{k: \sum_{r=1}^ky_r\geqslant n\}\}$.
\end{theorem}
\proof
Similar to the proof of part (ii) of Theorem \ref{circ}, let the vertex set of  $G$ be $\{v_1,\ldots , v_{\vert V(G)\vert}\}$, the vertex set of $ i$-th copy of $H$ be $\{w^{(i)}_1,\ldots , w^{(i)}_{\vert V(H)\vert}\}$ and  $ e_{ik} $ be the edge from $v_i$ to $w_k^{(i)}$. 
   We present an edge labeling of $G\circ H$ that is continuation of used edge labeling in the proof of part (ii) of Theorem \ref{circ}. We have the following steps:

{\bf Step 1)} We assign the edges  $e_{11},\ldots , e_{1m}$ the label $1$.  Set $x'_1=1$ and $y'_1=1$.

{\bf Step 2)} We assign the edges $e_{i1},\ldots , e_{i(i-1)}$ the label $2$ and $e_{ii},\ldots , e_{im}$ the label $1$, for $2\leqslant i \leqslant m$. Set $x'_2=m-1$.

{\bf Step 3)} Label The edges $e_{(m+1)1},\ldots , e_{(m+1)m}$ with the label $2$.

\medskip
So we used the label $2$ for labeling the edges between $G$ and $m$ copies of $H$. Set $y'_2=m$.
\medskip

{\bf Step 4)} Label $e_{(m+2)1},\ldots , e_{(m+2)m}$  with the label $3$. Next we do the same work as in Step $2$ with the label $1,3$ and $2,3$. So we labeled the edges $e_{i1},\ldots , e_{im}$ for $m+3\leqslant i \leqslant 3m$.

{\bf Step 5)} In this step we use the labels $1,2,3$ for the labeling of middle edges. We assign the first three edges $e_{i1},e_{i2},e_{i3}$ the labels $1,2,3$ for $3m+1\leqslant i \leqslant 3m+ x'_3$, where $x'_3= \sum_{j=2}^m (m-j)$. For labeling $e_{i4},\cdots , e_{im}$ we use the label $1,2,3$ such that $(L^{(i)}_1,L^{(i)}_2,L^{(i)}_3)$ are distinct for each $3m+1\leqslant i \leqslant 3m+ x'_3$, where $L^{(i)}_j$ is the number of the label $j$ in edges $e_{i1},\cdots , e_{im}$.  It can be seen that this number is $x'_3= \sum_{j=2}^m (m-j)$. 

\medskip
So we used the label $3$ for labeling the edges between $G$ and $1+2(m-1)+x'_3$ copies of $H$. Set $y'_3=1+2(m-1)+x'_3$.

By  continuing this method we get:
\begin{equation*}
x'_r=\left\{
\begin{array}{ll}
 1&r=1\\
 m-1&r=2\\
 \sum_{i_{r-2}=r-1}^{m}\cdots \sum_{i_2=i_3}^{m}\sum_{i_1=i_2}^{m}(m-i_1)&r\geqslant 3,
 \end{array}\right. 
\end{equation*}
\begin{equation*}
 y'_r=\left\{
 \begin{array}{ll}
 1&r=1\\
 m&r=2\\
 \sum_{i=0}^{r-1}{r-1\choose i}x_{i+1}&r\geqslant 3.
\end{array}\right.
\end{equation*}

With this method we have labeled (distinguishing)  all edges between the  vertices of $G$ and the vertices of copies of $H$. We use the label $1$ for the rest of edges. Therefore $D'(G\circ H)\leqslant min\{k: \sum_{r=1}^ky_r\geqslant n\}$.

On the other hand if we label $G$ in a distinguishing way with $D'(G)$ labels and assign the remaining edges the label $1$, then we obtain a distinguishing labeling of $G\circ H$ with $D'(G)$ labels, because $D'(H)=1$ and $H\neq K_2$. Therefore by the above paragraph we have the result. \qed

\medskip

\noindent {\bf Acknowledgement. } The authors would like to express their gratitude to the referee for her/his careful reading and
helpful comments which improved the paper.

\end{document}